\documentclass{article}

\usepackage{arxiv}

\usepackage[utf8]{inputenc} 
\usepackage[T1]{fontenc}    
\usepackage{hyperref}       
\usepackage{url}            
\usepackage{booktabs}       
\usepackage{amsfonts}       
\usepackage{nicefrac}       
\usepackage{microtype}      
\usepackage{lipsum}
\usepackage{graphicx}
\usepackage{amsmath}
\usepackage{caption}
\usepackage{subcaption}
\usepackage{amssymb}
\usepackage{algorithm}
\usepackage{tcolorbox} 
\graphicspath{ {./images/} }

\newtheorem{problem}{Problem}

\newtheorem{theorem}{Theorem}

\newtheorem{algo}{Algorithm}
\newtheorem{conclusion}{Conclusion}

\title{Time-Optimal Guidance for Intercepting Moving Targets with Impact-Angle Constraints}

\author{
 Yuan Zheng \\
  School of Aeronautics and Astronautics\\
  Zhejiang University\\
  Hangzhou 310027, Zhejiang, China\\
  \And
  Zheng Chen* \\
   School of Aeronautics and Astronautics\\
  Zhejiang University\\
  Hangzhou 310027, Zhejiang, China\\
  Corresponding author, Tel. +86-571-87953045\\
  \texttt{z-chen@zju.edu.cn} \\

}

\begin{document}
\maketitle
\begin{abstract}
The minimum-time path for intercepting a moving target with a prescribed impact angle is studied in the paper. The candidate paths from Pontryagin's maximum principle are analyzed, so that each candidate is related to a zero of a real-valued function. It is found that the real-valued functions or their first-order derivatives can be converted to polynomials of at most fourth degree. As a result, each canidate path can be computed within a constant time by embedding a standard polynomial solver into the typical bisection method. The control strategy along the shortest candidate eventually gives rise to the time-optimal guidance law. Finally, the developments of the paper is illustrated and verified by three numerical examples. 

\end{abstract}

\keywords{Dubins vehicle\and Minimum-time path\and Path planning\and Intercept guidance}

\section{Introduction}

 In this paper,  we consider a 2-dimensional pursuer-target engagement. The pursuer moves forward at a constant speed with a minimum turning radius, and the target moves at a constant speed. Such an engagement  is probably one of the most popular problems in the field of endgame guidance. In practical scenarios,  minimizing the engagement duration is crucial for the pursuer to successfully intercept the target, as it is essential to reduce probability of detection and  to improve  survivability against countermeasures \cite{Gopalan:2017}.  For this reason,  the Time-Optimal  Guidance Problem (TOGP) for a pursuer to intercept a target has  been widely studied in the literature. 

 It should be noted that the previously described pursuer takes the same kinematic model as the Dubins vehicle \cite{Dubins:57}. Therefore, when the target is stationary and the final impact angle is not constrained,  the TOGP is degenerate to the well-known Relaxed Dubins Problem  (RDP) \cite{Bui:94}. The solution path of RDP has been proven in \cite{Bui:94} to be in a sufficiently family of four candidates.  If  the final impact angle is fixed, it is known according to \cite{Dubins:57} that the time-optimal path can be computed within a constant time by checking at most six candidate paths.  


From practical point of view,  it is important to consider that the target is moving.  Without constraints on the final impact angle,  Mayer, Isaiah, and Shima \cite{Mayer:2015} established some sufficient conditions to ensure that the time-optimal path for intercepting a moving target shared the same  geometric pattern as the path of RDP. Whereas,  there was a gap between necessary and sufficient conditions, and it was not clear how to devise the time-optimal guidance law if the sufficient conditions were not met.  More recently,  the minimum-time paths for intercepting  moving targets were  thoroughly synthesized in  \cite{Yuan:2020}, which allowed developing  an efficient and robust algorithm  to compute the corresponding time-optimal  guidance law for intercepting moving targets.

If the target moves along a straight line, it can be proven by a simple coordinate transformation that the TOGP is equivalent to the problem of planning minimum-time path in a constant drift field \cite{Bakolas:2013}.  Many practical applications require addressing such path planning problems since the motions of  aerial vehicles and underwater vehicles are usually affected by wind and ocean current, respectively. Up to now, planning minimum-time paths in constant drift  has received some attentions.  For instance, without a constraint on final heading angle, the minimum-time paths in a constant drift field were related to zeros of some nonlinear equations, and the typical Newton iterative method and bisection method were proposed to find the minimum-time paths \cite{McGee:2007,McGee:2005,Techy:2009}.

A specific terminal heading is essential for various pursuer-target engagements. For example, a pursuer with directed warhead, against ground and ocean targets, is more effective when the impact occurs at a certain angle \cite{ Gopalan:2017}.  For this reason, the study on optimal fixed-impact-angle guidance is quite active in the field of endgame guidance. The TOGP for intercepting moving targets with lateral impact angle was first studied in \cite{Gopalan:2016}; it was found that the time-optimal guidance law was actually determined by a zero of a  highly nonlinear equation.  A natural extension was presented in \cite{Gopalan:2017} where the final impact angle could be assigned to any value. To be specific, by a coordinate transformation, the authors of \cite{Gopalan:2017}  managed to convert the TOGP with a general impact angle to the special problem in \cite{Gopalan:2016}. As a result, the method presented in \cite{Gopalan:2016} was used to find the time-optimal fixed-impact-angle guidance law by finding a zero of  nonlinear equations.


It is worth noting  that a nonlinear equation may have multiple zeros but only a specific zero is related to the time-optimal fixed-impact-angle guidance law.  Thus, the typical Newton-like iterative method and bisection method proposed in \cite{Gopalan:2016,Gopalan:2017,McGee:2007,McGee:2005} are not robust to find the optimal guidance law. The reasons include that 1) the numerical methods may not be able to converge to a zero if the initial guess is not appropriately chosen, and 2) even if the numerical methods converge to a zero, it is not necessarily the desired one related to the optimal guidance law, as shown by the numerical examples in Section \ref{SE:Numerical}. 
In this paper,  the solution path of the TOGP for intercepting a moving target with fixed impact angle is synthesized and some geometric properties are presented. Using these geometric properties, some  nonlinear equations in terms of the solution path's parameters are formulated. These nonlinear equations or their derivatives can be transformed to some polynomials of at most 4th degree. As a consequence, simply embedding a standard polynomial solver into the bisection method  leads to a robust and efficient method for finding the time-optimal fixed-impact-angle guidance law.  Since the problem of planning shortest path in constant drift field can be converted to the TOGP with a moving target, it follows that the developments of this paper also allow efficiently planning shortest paths in constant drift field. 
 
This paper is organized as follows. In Section \ref{SE:preliminary}, the TOGP for intercepting moving targets is formulated.  Necessary conditions for optimality are established and transcendental equations in terms of  the solution path's parameters are formulated in Section \ref{SE:Characterization}.  In Section 
\ref{SE:Analytical},  a numerical method is presented so that the time-optimal fixed-impact-angle guidance law can be computed within a constant time. Numerical examples are presented in Section \ref{SE:Numerical},  verifying and illustrating the developments of this paper. 
\section{Problem Formulation}\label{SE:preliminary}



Consider the 2-dimensional pursuer-target engagement scenario  presented in Fig.~\ref{Fig:geometry}. The inertial frame $Oxy$ is located in the horizontal plane, and the origin is the same as the initial position of the pursuer. The positive $x$-axis points to the east, and the $y$-axis is aligned with the north.  
\begin{figure}[ht]
	\centering
	\includegraphics[width = 2.5in]{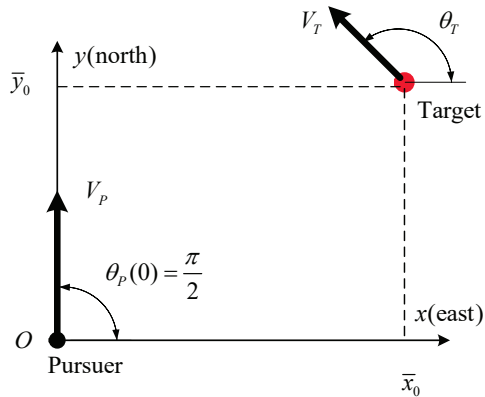}
	\caption{Geometry and coordinates system.}
	\label{Fig:geometry}
\end{figure}
The pursuer moves only forward at a constant speed $V_P>0$ with a bounded turning radius $\rho>0$. The heading angle $\theta_P\in[0,2\pi]$ defines the direction of the velocity of the pursuer, measured from east in a counter-clockwise direction. Denote by $(x_P,y_P)\in \mathbb{R}^2$ the position of the pursuer. Then, the motion of the pursuer is governed by 
\begin{align}
\begin{split}
\dot{x}_P &= V_P \cos \theta_P(t)\\
\dot{y}_P &= V_P \sin \theta_P(t)\\
\dot{\theta}_p &=  V_P \frac{u(t)}{\rho}
\end{split}
\label{Eq:problem1}
\end{align}
where $t\geq 0$ denotes time, the dot denotes the differentiation with respect to time,  and $ u\in [-1,1]$ is the control input representing the lateral acceleration of the pursuer.

The target is considered to be moving but nonmaneuvering. The constant speed is denoted by $V_T>0$. Denote by $\theta_T\in [0,2\pi]$ the heading angle of the target. It is apparent that $\theta_T$ remains constant throughout the engagement as the target does not maneuver. 
Let the position of the target at initial time be $(\bar{x}_{0},\bar{y}_{0})$. Then,  the position of the target at any time $t\geq 0$ is given by
\begin{align}
[x_T(t),y_T(t)] = [\bar{x}_{0},\bar{y}_{0}]+ V_Tt[\cos\theta_T,\sin\theta_T]
\nonumber
\end{align}
Without loss of generality, we assume that the initial heading angle of the pursuer is $\theta_{P_0}=\pi/2$ so that the state of the pursuer at initial time $t=0$ is
$$\boldsymbol{z}_0 := (0,0,\frac{\pi}{2})$$
 Then, finding the time-optimal fixed-impact-angle guidance law for intercepting a moving target is equivalent to addressing the following Optimal Control Problem (OCP).
\begin{problem}[OCP]\label{problem1}
Find a minimum time $t_f > 0$ so that  the system in Eq.~(\ref{Eq:problem1}) is steered by a measurable control $u(\cdot)\in [-1,1]$ over the interval $[0,t_f]$ from the fixed initial state $\boldsymbol{z}_0$ at $t=0$ to intercept the target at $t_f$ with a fixed impact angle $\phi_f\in [0,2\pi)$, i.e., \begin{align}
\begin{split}
[x_P(t_f),y_P(t_f)] &= [x_T(t_f),y_T(t_f)]\\
\phi_f & = \theta_P(t_f) - \theta_T
\end{split}
\label{EQ:terminal_boundary}
\end{align}
\end{problem}
According to Eq.~(\ref{EQ:terminal_boundary}),  given the terminal impact angle $\phi_f$, the final heading angle of the pursuer is fixed as well, i.e.,
$$\theta_{P_f} := \phi_f + \theta_T$$ 
Throughout the paper, we assume that the speed ratio $V_T/V_P$ is less than $1$, which ensures that the solution of the OCP exists. \cite{Cockayne:1967}




\section{Characterizing the Solution of the OCP}\label{SE:Characterization}

Let  $\lambda_x,\ \lambda_y$ and $\lambda_{\theta}$ be the costate variables of $x_P,\ y_P$, and $\theta_P$, respectively. Then, the Hamiltonian of the OCP  is
\begin{align}
H=\lambda_x V_P \cos \theta_P +\lambda_y  V_P\sin \theta_P +\lambda_{\theta} V_P u/\rho-1
\label{Eq:problem4}
\end{align}
According to Pontryagin's maximum principle \cite{Pontryagin}, we have
\begin{align}
\dot{\lambda}_x(t)&=-\frac{\partial H}{\partial x_P}=0\label{EQ:px}\\
\dot{\lambda}_y(t)&=-\frac{\partial H}{\partial y_P}=0\label{EQ:py}\\
\dot{\lambda}_{\theta}(t)&=-\frac{\partial H}{\partial \theta_P}=\lambda_x(t) V_P\sin \theta_P(t)-\lambda_y(t)V_P \cos \theta_P(t)
\label{Eq:problem5}
\end{align}
It is apparent from Eq.~\eqref{EQ:px} and Eq.~\eqref{EQ:py} that $\lambda_x$ and $\lambda_y$ are constant. By integrating  Eq.~(\ref{Eq:problem5}), we have
\begin{align}
\lambda_{\theta}(t)=\lambda_x y_P(t)- \lambda_y x_P(t)+c_0
\label{Eq:problem6}
\end{align}
where $c_0$ is a scalar constant.  In view of Eq.~(\ref{Eq:problem6}),  if $\lambda_{\theta}\equiv 0$ on a nonzero interval, the path $(x,y)$ is a straight line segment on this interval. Note that $u\equiv 0$ along a straight line.  Thus,  we have  $u\equiv0$ if $\lambda_{\theta} \equiv 0$.  As a result, the maximum principle indicates that the optimal control $u$ is totally determined by $\lambda_{\theta}$, i.e.,
\begin{align}
u=\begin{cases}
1,&\lambda_{\theta}>0\\
0,&\lambda_{\theta}\equiv0\\
-1,&\lambda_{\theta}<0
\end{cases}
\label{Eq:problem9}
\end{align}
The path $[x(t),y(t)]$ is a circular arc with right (resp. left) turning direction if $u = -1$ (resp. $u=1$). Therefore, the switching conditions in Eq.~\eqref{Eq:problem9} imply that  the solution path of the OCP is a concatenation of circular arcs and straight line segments.

It has been proven in \cite{Dubins:57} that, if the distance between the initial and final positions of the pursuer is at least $4\rho$, the solution of the OCP is a circular arc, followed by a straight line segment, and followed by another circular arc. In this paper, the separation between the initial and final positions  is also assumed to be at least $4\rho$, as was done in the seminal works \cite{Gopalan:2016,Gopalan:2017}. Then, the geometric pattern of the OCP's solution can be denoted by CSC, where ``C'' and ``S'' represent a circular arc with radius of $\rho$ and  a straight line segment, respectively.   If a circular arc C has a right (resp. left) turning direction, we represent it by R (resp. L). Then, we have that the CSC type includes four different types, i.e., $$\mathrm{CSC} = \{\mathrm{RSR},\ \mathrm{RSL},\ \mathrm{LSR},\ \mathrm{LSL}\}$$


Note that the speed of pursuer is constant, and the value of control on each subarc is available. Therefore, in order to find the solution path of the OCP, it amounts to finding the randian of each circular arc and the length of the straight line segment.  For notational simplicity,  we use the notation C$_{\alpha}$ to denote a circular arc of radian $\alpha\geq 0$, and use the notation S$_{d}^{\beta}$ to denote a straight line segment of length $d\geq 0$ where $\beta\in [0,2\pi)$ denotes the orientation angle of the straight line segment with respect to $x$-axis.  In the remainder of this paper, we use $C_{\alpha}S_{d}^{\beta}C_{\gamma}$ to represent CSC when necessary.  All the four types in $C_{\alpha}S_{d}^{\beta}C_{\gamma}$ are illustrated in Fig.~\ref{Fig:geometry}.



\begin{figure}[!htp]
\centering
\begin{subfigure}[t]{6cm}
\centering
\includegraphics[width = 6cm]{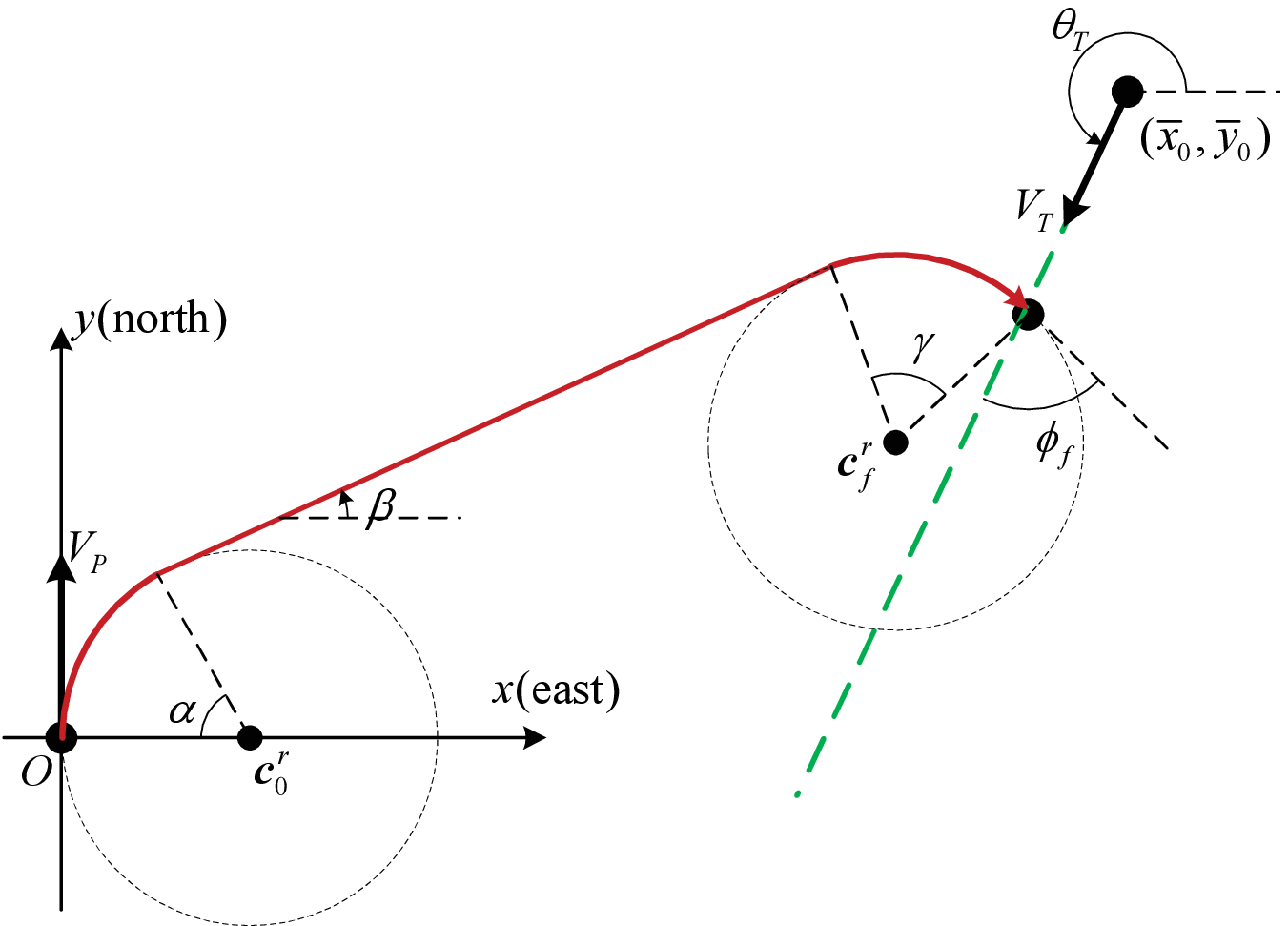}
\caption{RSR}
\label{Fig:geometry_RSR}
\end{subfigure}
~~~~~
\begin{subfigure}[t]{6cm}
\centering
\includegraphics[width = 6cm]{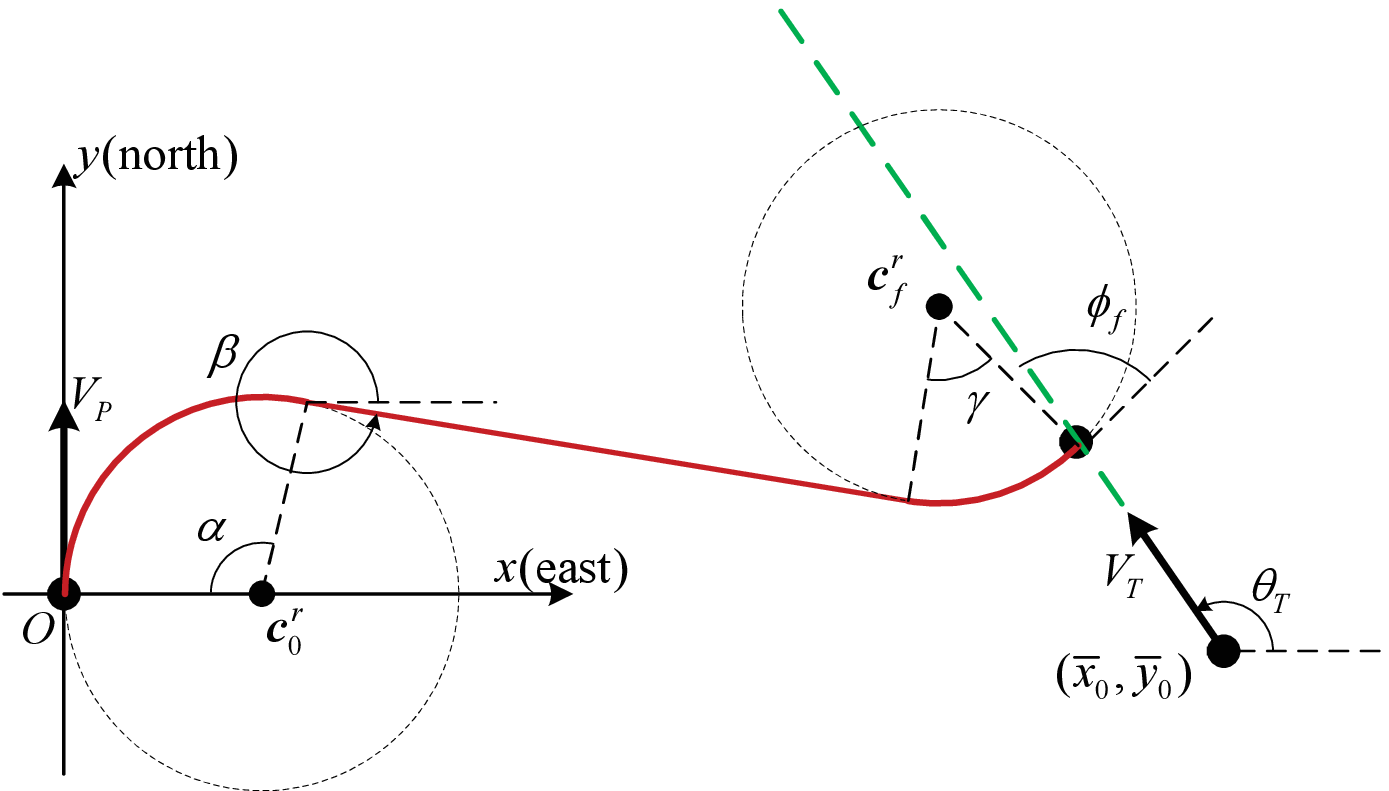}
\caption{RSL}
\label{Fig:geometry_RSL}
\end{subfigure}\\
\begin{subfigure}[t]{6cm}
\centering
\includegraphics[width = 6cm]{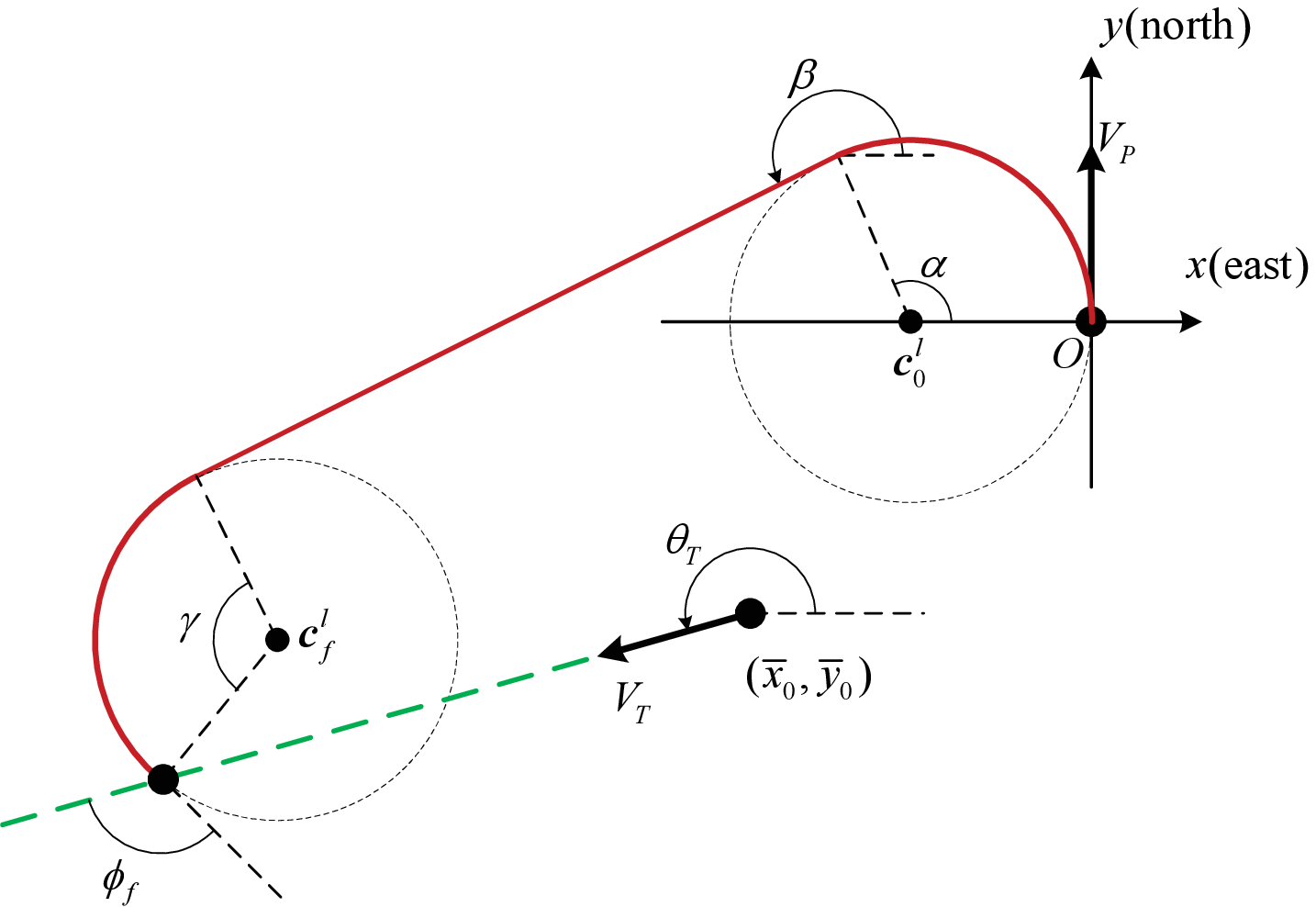}
\caption{LSL}
\label{Fig:geometry_LSL}
\end{subfigure}
~~~~~
\begin{subfigure}[t]{6cm}
\centering
\includegraphics[width = 6cm]{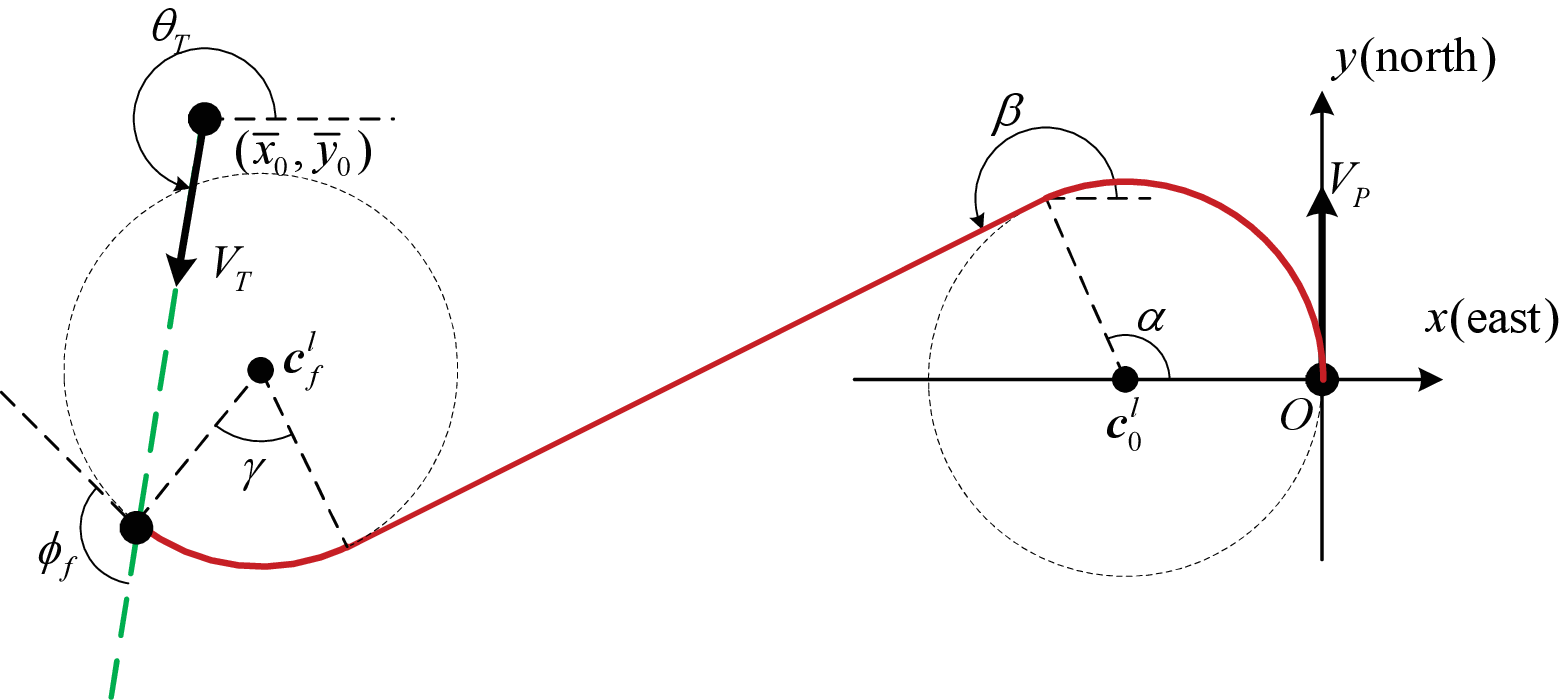}
\caption{LSR}
\label{Fig:geometry_LSR}
\end{subfigure}
\caption{Geometry for the paths of C$_{\alpha}$S$_{d}^{\beta}$C$_{\gamma}$.}
\label{Fig:geometry}
\end{figure}



By simple geometric analysis, we  have that the values of $\alpha$, $d$, and $\gamma$ are totally determined by $\beta$, as shown in Appendix \ref{Appendix:A}. Therefore, once the value of $\beta$ is obtained, we immediately have the solution path of the OCP, which gives rise to the optimal control strategy or the time-optimal fixed-impact-angle guidance law. 
By the following theorem, it is shown that $\beta$ is a zero of some nonlinear equations.

\begin{theorem}\label{TH:1}
If the solution path of the OCP is of type C$_{\alpha}$S$_d^{\beta}$C$_{\gamma}$,  the following four statements hold:
\begin{description}
\item (1) If both $C_{\alpha}$ and $C_{\gamma}$ are right-turning circular arcs, we have
\begin{align}
a_1 + a_2 \sin \beta +a_3\cos \beta  =0
\label{EQ:RSR}
\end{align}
where $a_1$--$a_3$ are constants  given in Appendix \ref{Appendix:B}.
\item (2) If both $C_{\alpha}$ and $C_{\gamma}$ are left-turning circular arcs, we have
\begin{align}
b_1 +b_2 \sin \beta + b_3\cos \beta  =0
\label{EQ:LSL}
\end{align}
where $b_1$--$b_3$ are constants given in Appendix \ref{Appendix:B}.
\item (3) If $C_{\alpha}$ is a right-turning circular arc and $C_{\gamma}$ is a left-turning circular arc, we have
\begin{align}\label{EQ:RSL}
c_1 + c_2 \sin \beta +c_3 \cos \beta + \beta(c_4 \sin \beta +c_5 \cos \beta) = 0
\end{align}
where $c_1$--$c_5$ are constants given in Appendix \ref{Appendix:B}.
\item (4) If $C_{\alpha}$ is a left-turning circular arc and $C_{\gamma}$ is a right-turning circular arc, we have
\begin{align}\label{EQ:LSR}
d_1 + d_2 \sin \beta +d_3 \cos \beta + \beta(d_4 \sin \beta +d_5 \cos \beta) = 0
\end{align}
where $d_1$--$d_5$ are constants given in Appendix \ref{Appendix:B}.
\end{description}
\end{theorem}
The proof of this theorem is postponed to Appendix \ref{Appendix:B}.  In the next section, the results in Theorem \ref{TH:1} will be used to find the solution of the OCP. 

\section{Numerical Method for Finding the Optimal Control}\label{SE:Analytical}

As presented in Appendix \ref{Appendix:A}, it suffices to compute the value of $\beta$ in order to devise the optimal control strategy.  According to Theorem \ref{TH:1},  if the solution path is of type RSR or LSL,  the value of $\beta$ can be obtained analytically by finding the zeros of Eq.~(\ref{EQ:RSR}) and Eq.~(\ref{EQ:LSL}), respectively.  However, the transcendental equations Eq.~(\ref{EQ:RSL}) and Eq.~(\ref{EQ:LSR}) may have multiple zeros but only a specific zero is related to the optimal path. Existing numerical solvers cannot be guaranteed to find the desired zero related to the optimal path, as illustrated by the numerical examples in Section \ref{SE:Numerical}.   In the remainder of this section, a robust and efficient method will be presented to find all the zeros of equations having the same form as Eq.~(\ref{EQ:RSL}) and Eq.~(\ref{EQ:LSR}). As a result, the optimal path of the OCP can be computed robustly and efficiently by ruling out useless zeros.

Let us define a function 
\begin{align}
G(\beta)  \overset{\triangle}{=}  e_1 \sin \beta + e_2 \cos \beta + \beta(e_3 \cos \beta+e_4\sin \beta) + e_5
\label{EQ:new_fun}
\end{align}
where  $e_1$--$e_5$ are scalar constants.  It is apparent that both Eq.~(\ref{EQ:RSL}) and Eq.~(\ref{EQ:LSR}) have the same form as Eq.~(\ref{EQ:new_fun}).   

Set 
\begin{align}
G_1(\beta)  & \overset{\triangle}{=} \beta + ({e_5 + e_1 \sin \beta + e_2 \cos \beta})/({e_4\sin \beta + e_3 \cos \beta})\nonumber\\
G_2(\beta) &\overset{\triangle}{=} e_5 + e_1 \sin \beta + e_2 \cos \beta\nonumber
\end{align}
 Then, by rearranging Eq.~(\ref{EQ:new_fun}), we have that the zeros of $G(\beta)$ are equivalent to those of
\begin{align}\label{EQ:trans}
\bar{G}(\beta) \overset{\triangle}{=}
\begin{cases}
G_1(\beta), \ \text{if}\ e_4\sin \beta + e_3 \cos \beta \neq 0\\
G_2(\beta),\  \text{if}\ e_4\sin \beta + e_3 \cos \beta = 0
\end{cases}
\end{align}
Since the zeros of $G_2(\beta)$ is readily available, the following paragraph will only be contributed to finding the zeros of $G_1(\beta)$.  

Differentiating $G_1(\beta)$ with respect to $\beta$ leads to
\begin{align}\label{EQ:pre_poly}
G_1^{\prime}(\beta)= &\ 1+   \big[(e_1 \cos \beta  - e_2 \sin \beta )(e_3 \cos \beta
 +   e_4 \sin \beta)   -   (e_1 \sin \beta + e_2 \cos \beta + e_5) \\
 \times &\  ( e_4 \cos \beta - e_3 \sin \beta )\big]/(e_3 \cos \beta
  +\  e_4 \sin \beta)^2\notag
\end{align}
By substituting the half-angle formulas 
\begin{align}
\sin \beta = \frac{2\tan \frac{\beta}{2}}{1+ \tan^2\frac{\beta}{2}} \ \text{and}\ \cos \beta = \frac{ 1- \tan^2 \frac{\beta}{2}}{1+ \tan^2\frac{\beta}{2}}
\end{align}
into Eq.~(\ref{EQ:pre_poly}), we have that the zeros of $G_1^{\prime}(\beta)$ are equivalent to those of the following quartic polynomial:
\begin{align}\label{EQ:polynomial}
p_1 \tan^4(\frac{\beta}{2}) + p_2 \tan^3(\frac{\beta}{2}) + p_3 \tan^2(\frac{\beta}{2}) + p_4 \tan(\frac{\beta}{2}) + p_5 = 0
\end{align}
where
\begin{align}
\begin{split}
p_1 & =e_3^2+e_1e_3+e_4e_5-e_2e_4 \\
p_2 & = -4e_3e_4+2e_3e_5\\
p_3 &=2e_1e_3-2e_2e_4-2e_3^2+4e_4^2\\
p_4 & = 4e_3e_4+2e_3e_5\\
p_5 & =e_3^2+e_1e_3-e_2e_4-e_4e_5
\end{split}\notag
\end{align}
Since the roots of any quartic polynomial can be readily obtained either by radicals or by standard polynomial solvers, it follows that the zeros of $G_1^{\prime}(\alpha)$ can be obtained immediately. 

Note that the differentiation $G_1^{\prime}(\beta)$ has at most $4$ real zeros. Let us denote all the real zeros of $G^{\prime}_1(\beta)$ by $\beta_1$, $\ldots$, $\beta_n$  where $n\leq 4$.  Without loss of generality, we assume $\beta_0<\beta_1<  \ldots < \beta_n<\beta_{n+1}$ where $\beta_0 = 0$ and $\beta_{n+1}=2\pi$. Then, according to  \cite[Lemma 6]{Yuan:2020},  we have the following two conclusions.
\begin{conclusion}\label{conclusion1}
 For any $i\in \{0,1,\ldots,n\}$,  if $G_1(\beta_i)\times G_1(\beta_{i+1}) < 0$, the function $G_1(\beta)$ over the interval $[\beta_{i},\beta_{i+1}]$ has only one zero.
\end{conclusion}
\begin{conclusion}\label{conclusion2}
For any $i\in \{0,1,\ldots,n\}$,  if  $G_1(\beta_i)\times G_1(\beta_{i+1}) > 0$, the function $G_1(\beta)$ over the interval $[\beta_{i},\beta_{i+1}]$ does not have a zero. 
\end{conclusion}

According to Conclusion \ref{conclusion1}, if $G(\beta_i)\times G(\beta_{i+1})<0$, we can use a simple bisection method to find the only zero in the interval $(\beta_{i},\beta_{i+1})$. We denote by 
\begin{align}
z = B[G(\beta),\beta_i,\beta_{i+1}]
\end{align}
the bisection method to find the zero $z$ of $G(\beta)$ in the interval $(\beta_i,\beta_{i+1})$. With these notations, we can obtain all the real zeros of $G_1(\beta)$ in Eq.~(\ref{EQ:trans})  by the procedure in Algorithm \ref{algo1}.

\begin{center}
\begin{tcolorbox}[
colframe=blue!25,
colback=blue!10,
coltitle=blue!20!black,  
fonttitle=\bfseries,
width = 0.6\textwidth]
\begin{center}
\begin{algo}
[Finding all the zeros of $G_1(\beta)$]\label{algo1}
\end{algo}
\end{center}
\begin{description}
\item step 0. Set $i=0$ and $Z= \emptyset$. 
\item step 1. If $i\leq n$, go to step 2; otherwise, go to step 3.
\item step 2.  If $G_1(\beta_i) = 0$
\begin{description}
\item   \ \ \ \ \ \  \ \ \ \ \ \ \ \ \ \ $Z = Z\cup \{\beta_i\}$
\end{description}
\item \ \ \ \ \  \ \ \ \ \ \ Elseif\ $G_1(\beta_i)\times G(\beta_{i+1}) < 0$
\begin{description}
\item     \ \ \ \ \ \   \ \ \ \ \  \ \ \ \ \  $ z = B[G_1(\beta),\beta_i,\beta_{i+1}]$
\item   \ \ \ \ \ \   \ \ \ \ \  \ \ \ \ \  $Z = Z\cup \{z\}$
\end{description}
\item  \ \ \ \ \  \ \ \ \ \ \  Endif
\item  \ \ \ \ \  \ \ \ \ \ \  Set $i=i+1$ and go to step 1.
\item step 3. End
\end{description}
\end{tcolorbox}
\end{center}

Let us gather a few words to explain the pseudo codes in Algorithm \ref{algo1}. If $G_1(\beta_i) = 0$, we have that $\beta_i$ is a zero of $G_1(\beta)$; thus, we add $\beta_i$ into the set $Z$.  If $G_1(\beta_i)\times G_1(\beta_{i+1})<0$, in view of Conclusion \ref{conclusion1} we have that the function $G_1(\beta)$ on $(\beta_i,\beta_{i+1})$ has only one zero. Thus, a typical bisection method can be used to find the zero, and it is added to the set $Z$. If $G_1(\beta_i)\times G_1(\beta_{i+1})>0$, none zero exists between $\beta_i$ and $\beta_{i+1}$ according to Conclusion \ref{conclusion2}. Thus, nothing is done in step 2 if $G_1(\beta_i)\times G_1(\beta_{i+1})>0$. It should be noted that the set $Z$ contains all the real zeros of $G_1(\beta)$ on $[0,2\pi]$.  According to above analysis, it is clear that Algorithm \ref{algo1} is robust to find all the zeros of $G_1(\beta)$ within a constant time. 

Note that the zeros of $G_2(\beta)$ in Eq.~(\ref{EQ:trans}) can be found analytically.  Thus, by employing Algorithm \ref{algo1} to find the zeros of $G_1(\beta)$,  all the zeros of $\bar{G}(\beta)$   can be obtained.  Since the zeros of $G(\beta)$ are equivalent to those of $\bar{G}(\beta)$,  it follows that  all the zeros of $G(\beta)$ can be found.  


Up to now, we are able to find all the zeros of Eqs.~(\ref{EQ:RSR}--\ref{EQ:LSR}) either by analytical method or by Algorithm \ref{algo1}.  Each zero is related to a candidate path of the OCP according to Appendix \ref{Appendix:A}. Obviously,  the shortest candidate path gives rise to the optimal path of the OCP. 

\section{Numerical Examples}\label{SE:Numerical}
In this section, some numerical examples were simulated to demonstrate the developments of this paper.  Before proceeding, it is worth mentioning that a large number of the examples were tested on a desktop with AMD Ryzen 2500U, showing that the time-optimal fixed-impact-angle guidance law for any example could be established within $10^{-4}$ seconds. 

For numerical convenience,  the position was normalized so that the speed of pursuer was one, i.e., $V_P = 1$ m/s.  



\subsection{Case A: Time-Optimal Guidance for Intercepting a Moving Target}
For case A,   a numerical example in \cite{Gopalan:2016} was chosen to demonstrate the developments of the paper. In the normalized setting,  the target's initial position and velocity are equivalent to $(\bar{x}_0,\bar{y}_0)=(3.996,-5.388)$ m and $(-0.492,-0.0868)$ m/s, respectively.  The heading angle of the pursuer at the final impact time is  equivalent to $\theta_{P_f}=-1.396$ rad. 

\begin{figure}[!htp]
\centering
\begin{subfigure}[t]{0.48\textwidth}
\centering
\includegraphics[width = 1\textwidth,trim=2.5cm 0cm 0cm 0cm]{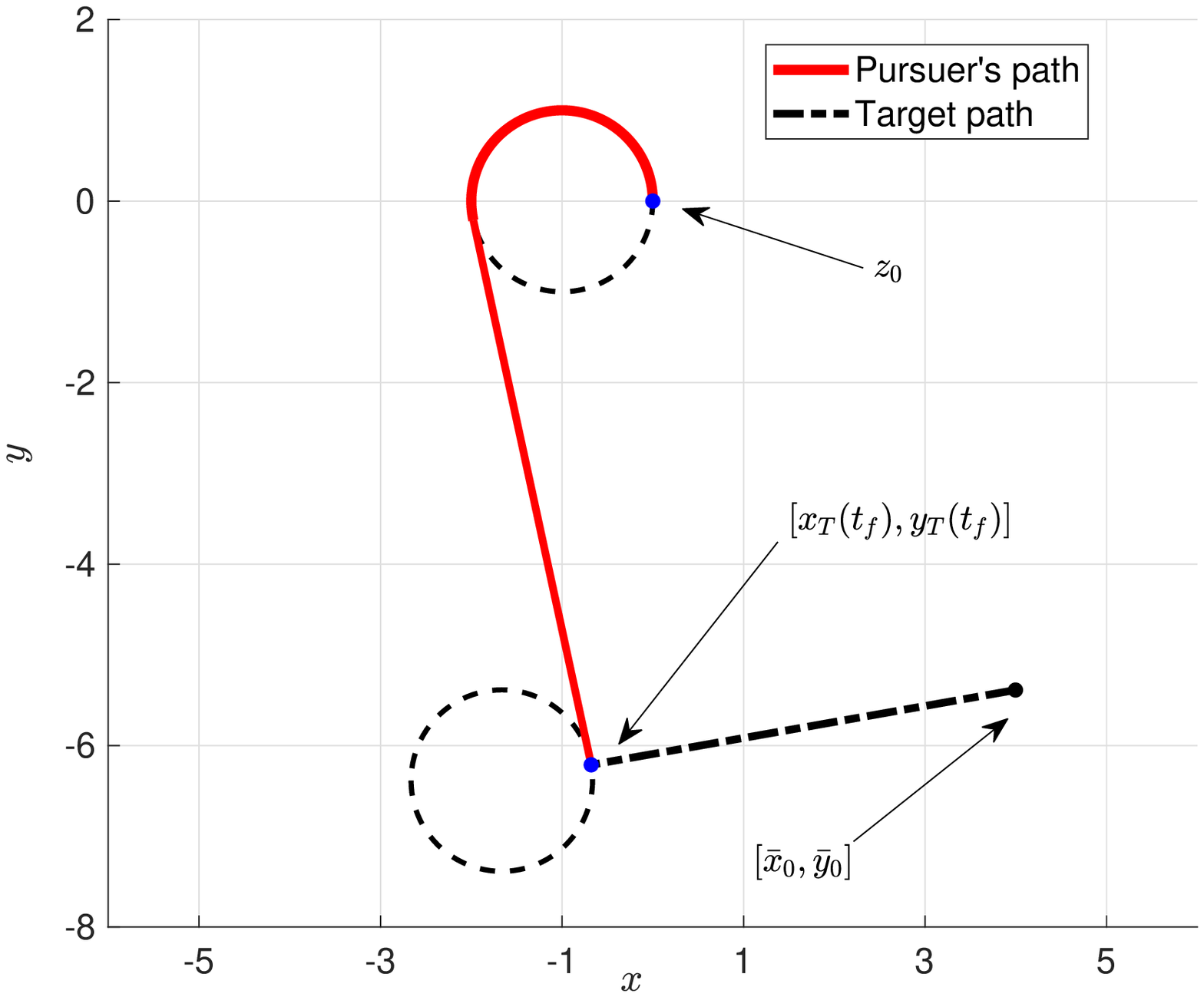}
\caption{Optimal path}
\label{Fig:case2a}
\end{subfigure}
\begin{subfigure}[t]{0.48\textwidth}
\centering
\includegraphics[width = 1\textwidth,trim=2cm 0cm 0cm 0cm]{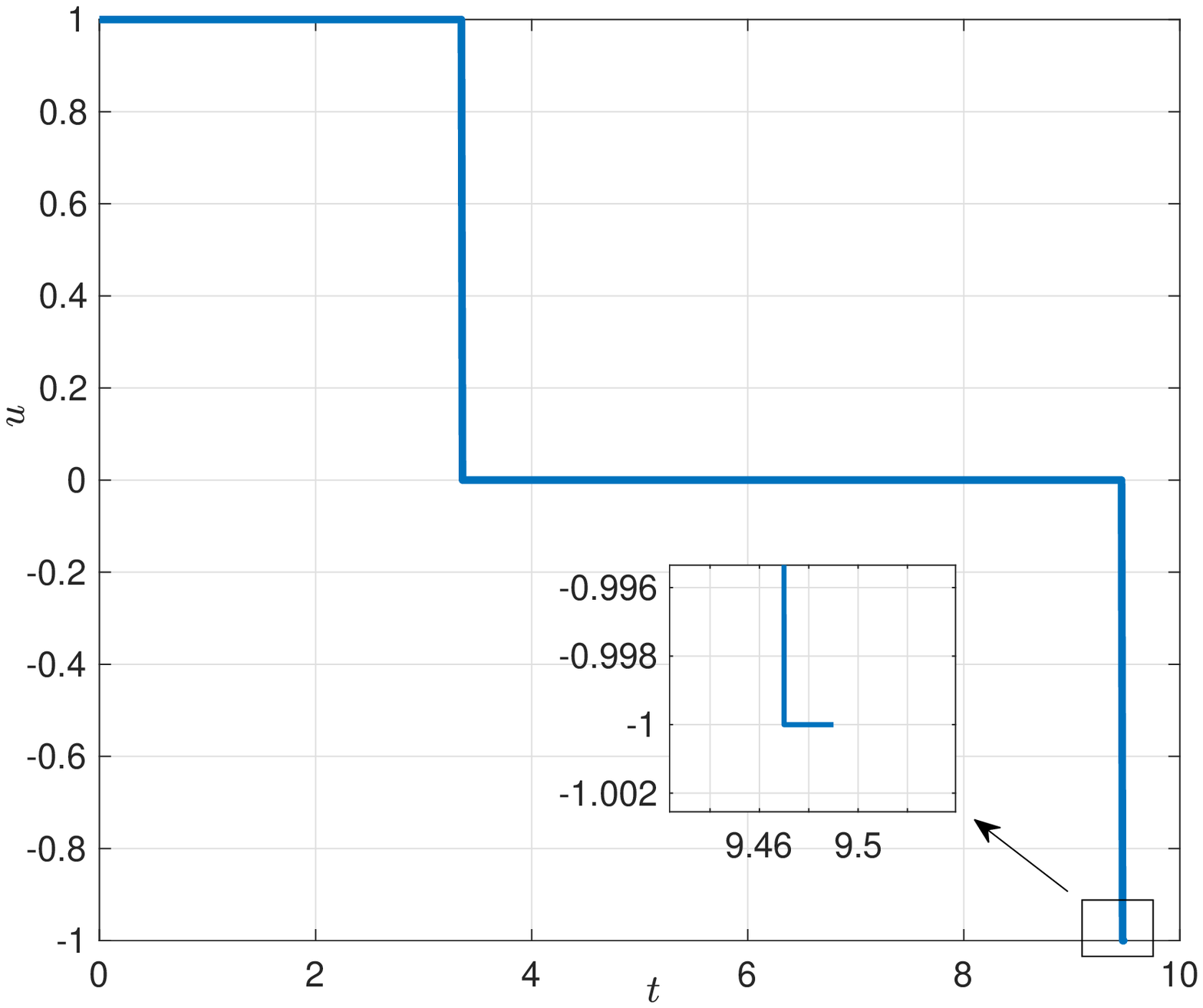}
\caption{Optimal control strategy}
\label{Fig:case2c}
\end{subfigure}
\caption{Case A: The solution path and corresponding optimal control strategy of the OCP.}
\label{Fig:caseB}
\end{figure}

In order to find the time-optimal path for case A,  one needs to check the lengths of all the four candidate paths in CSC, and the shortest candidate gives rise to the solution path.  Regarding the candidate path of RSR,  it was found analytically that Eq.~(\ref{EQ:RSR}) had two real zeros. Considering the fact that the two zeros represent the orientation angle of straight line segment, we can use a simple geometric analysis to rule out one useless zero, and another zero can be used to compute the length of path of RSR. For this example, the length of the candidate path of RSR was obtained as 16.63 m. Analogously,  by analytically finding the zeros of  Eq.~(\ref{EQ:LSL}), the length of the path of LSL was obtained as 9.58 m. When computing the lengths of the paths of RSL and LSR, Algorithm \ref{algo1} was used to find all the zeros of Eq.~(\ref{EQ:RSL}) and Eq.~(\ref{EQ:LSR}).  Numerical results indicated that both Eq.~(\ref{EQ:RSL}) and Eq.~(\ref{EQ:LSR}) had 2 real zeros. By geometric analysis, useless zeros were ruled out and the lengths of the paths of RSL and LSR were obtained as 9.49 m and 16.56, respectively.  Therefore, it is concluded that the shortest path for case A is of type RSL.  The shortest path is presented in
Fig.~\ref{Fig:case2a}, and the corresponding optimal control strategy is presented in Fig.~\ref{Fig:case2c}. The last circular subarc exists  but it is quite short, as shown by the scaled plot in Fig.~\ref{Fig:case2c}.

\begin{figure}[!hpt]
\centering
\includegraphics[width = 0.7\textwidth,trim=2cm 0cm 2cm 2cm]{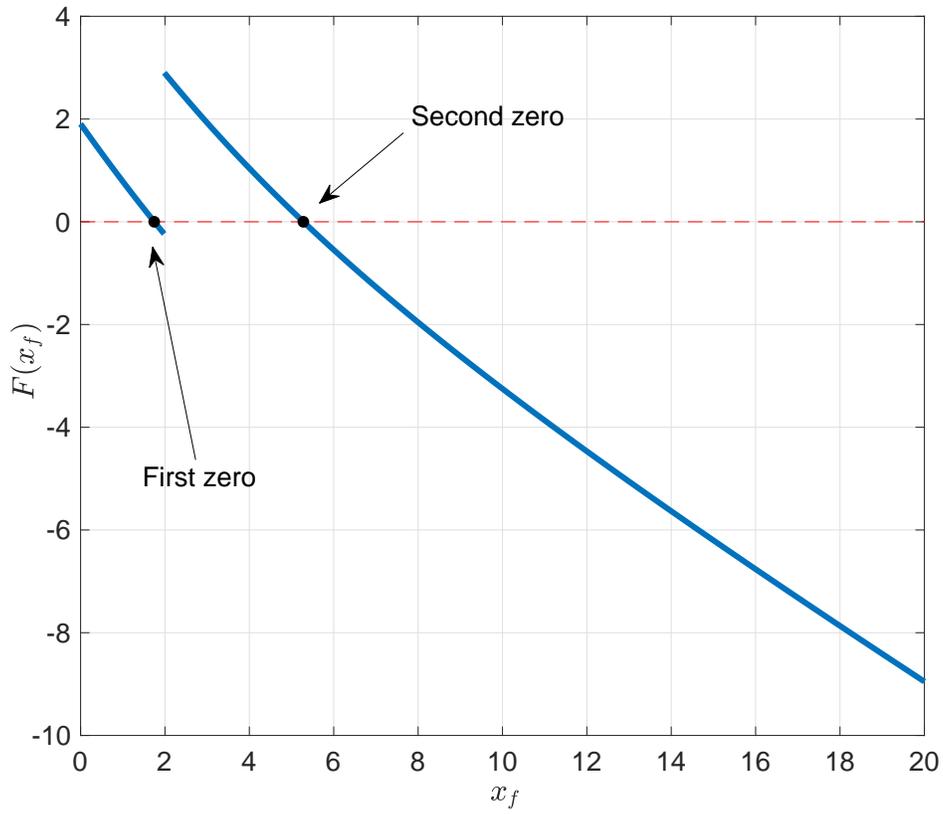}
\caption{CaseA: The function $F(x_f)$ related the path type of RSL in \cite{Gopalan:2016}.}
\label{Fig:case2b}
\end{figure}

\begin{figure}[!hpt]
\centering
\includegraphics[width = 0.6\textwidth,trim=2cm 0cm 2cm 1cm]{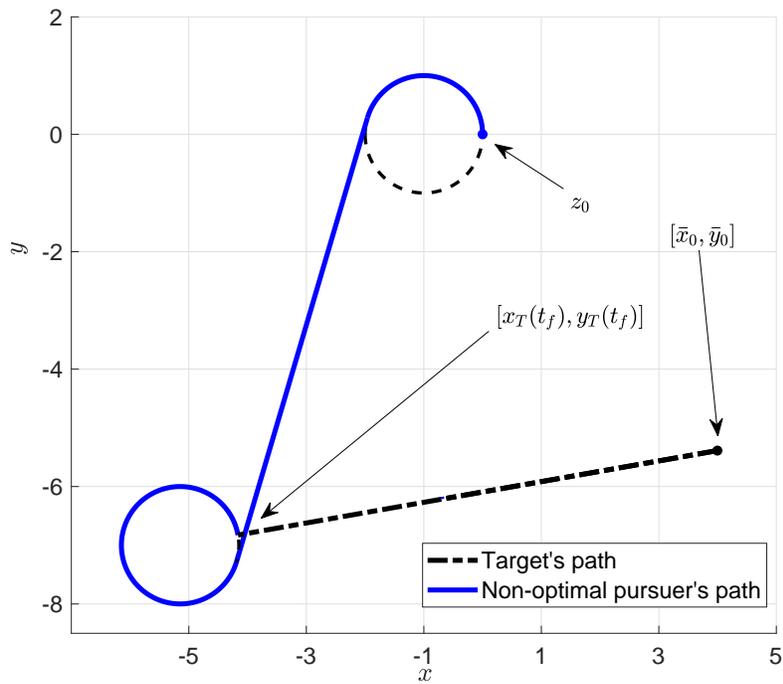}
\caption{Case A: Non-optimal path related to the second zero of $F(x_f)$ in \cite{Gopalan:2016}.}
\label{Fig:case2b1}
\end{figure}

Analogous to the present paper, it was proposed in \cite{Gopalan:2016} to find the optimal path  of this example by finding zeros of a nonlinear function $F(x_f)$ where $x_f$ is the projection of the final position on $x$-axis.  The value of the nonlinear function $F(x_f)$  against $x_f$ is ploted in Fig.~\ref{Fig:case2b}, showing that it is discontinuous and have more than one zero. 
Only is a specific zero related to the optimal path, but existing numerical solvers do not necessarily converge to the desired zero.  Once the second zero in Fig.~\ref{Fig:case2b} is found, it leads to a non-optimal path (the length is 16.56 m),  as presented in Fig.~\ref{Fig:case2b1}.



\subsection{Case B: Path Planning in Constant Drift Field}


The motions of aerial vehicles and underwater vehicles are usually affected by wind and ocean current, respectively. Thus, it is of practical importance to plan shortest path in a constant drift field. 


Take the path planning problem  in \cite{McGee:2005} to intercept a stationary target in a constant drift field as an example. In the normalized setting,  the constant velocity of the drift   is $\boldsymbol{w}=(0.3536, 0.3536)$ m/s.  The position of the stationary target is $(\bar{x}_0,\bar{y}_0)=(2.8284, 4.2426)$ m, and the final heading angle is $\theta_{P_f}=3.927$ rad. It was stated in \cite{McGee:2005} that this path planning problem could be solved by using Newton-like iterative method or bisection method to find a zero of some nonlinear equations. However, the formulas of nonlinear equations were not given in \cite{McGee:2005}. In addition, as stated above, even if the explicit expression of nonlinear equations are available, these numerical methods cannot be guaranteed to find a desired zero of nonlinear equations. In the next paragraph, we shall show how to find the shortest path in constant drift by the method proposed in this paper. 

\begin{figure}[!htp]
\centering
\begin{subfigure}[t]{0.48\textwidth}
\centering
\includegraphics[width = 1\textwidth,trim=2cm 0cm 2cm 0cm]{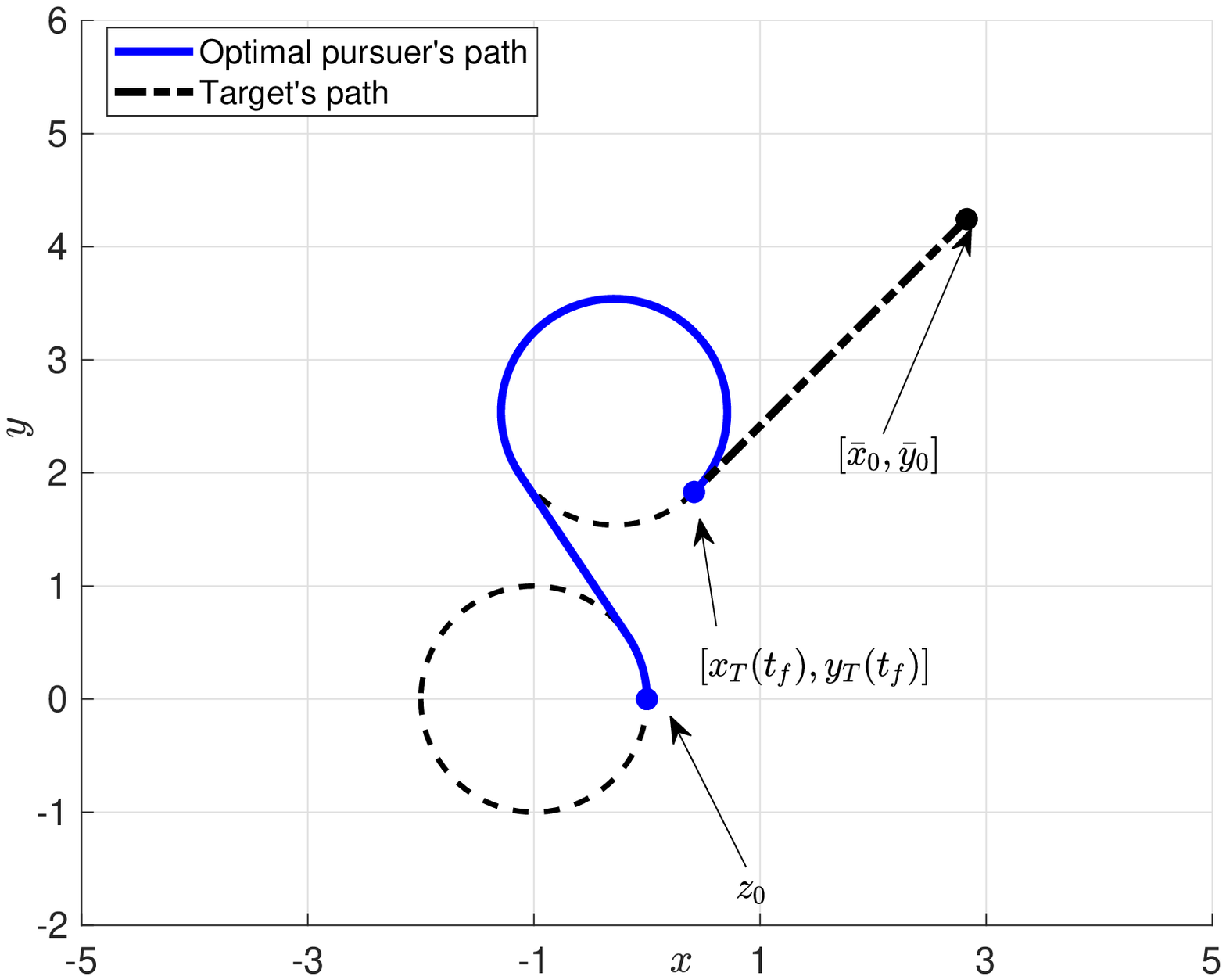}
\caption{Target's velocity is $-\boldsymbol{w}$ without constant drift}
\label{Fig:case3a}
\end{subfigure}
\begin{subfigure}[t]{0.48\textwidth}
\centering
\includegraphics[width =1\textwidth,trim=2cm 0cm 2cm 0cm]{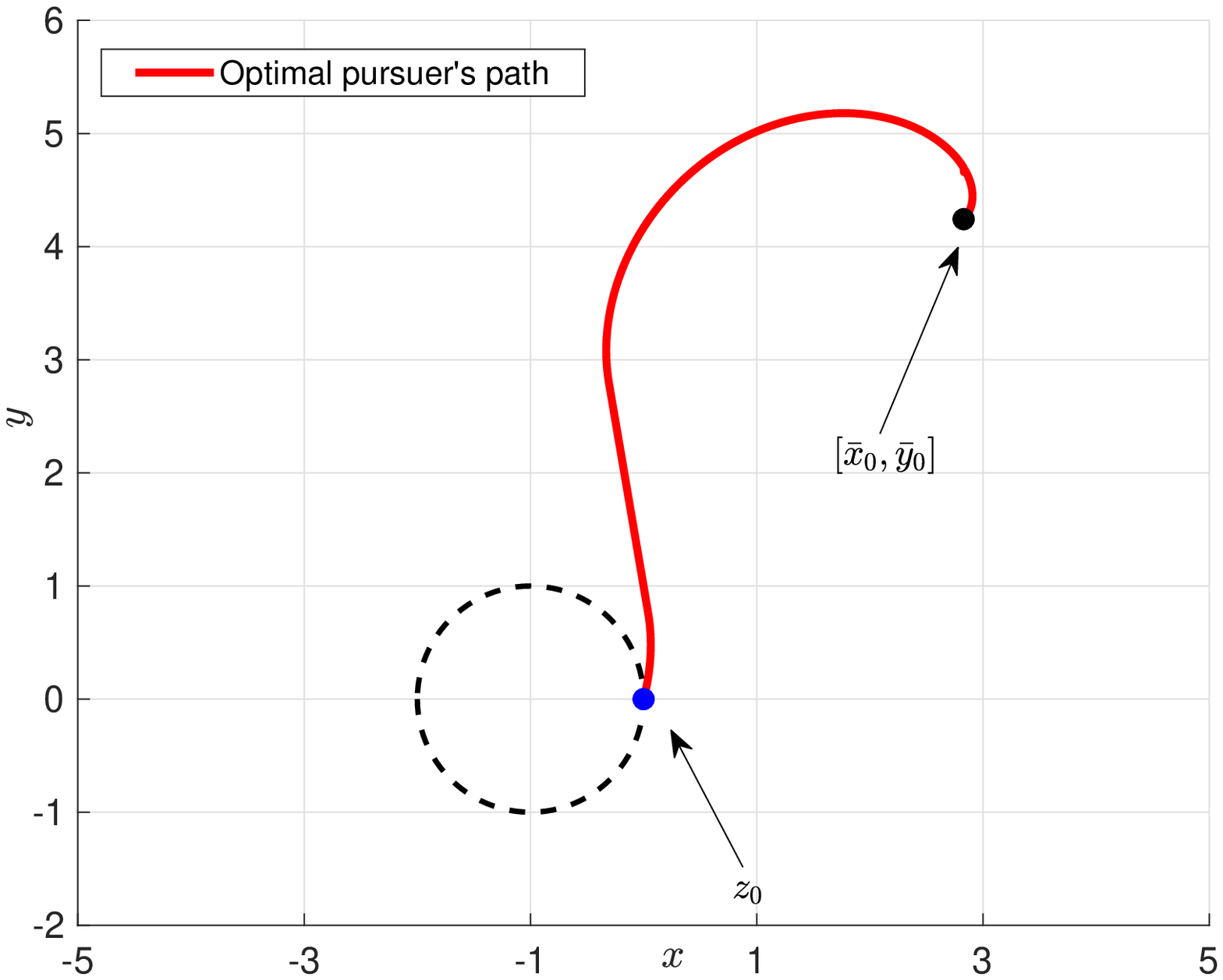}
\caption{Target's velocity is $\boldsymbol{0}$ with constant drift vector $\boldsymbol{w}$}
\label{Fig:case3a1}
\end{subfigure}
\caption{Case B: The solution path of the OCP and corresponding path in constant drift.}
\end{figure}
\begin{figure}[!htp]
\centering
\includegraphics[width = 0.6\textwidth,trim=2cm 0cm 2cm 0cm]{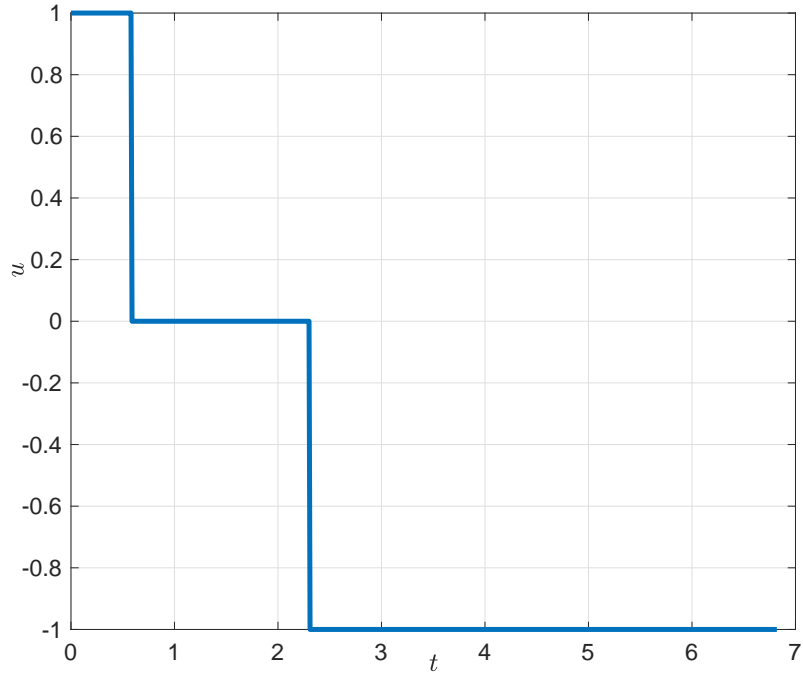}
\caption{Case B: The optimal control strategy against time.}
\label{Fig:case3c}
\end{figure}

By a simple coordinate transformation \cite{Yuan:2020,Bakolas:2013}, this path planning problem in a constant drift field can be equivalent to the TOGP with a moving target, and the velocity of the moving target is equivalent to $-\boldsymbol{w}$.  According to the numerical procedure in Section \ref{SE:Analytical}, by finding zeros of nonlinear equations in Eqs.~(\ref{EQ:RSR}--\ref{EQ:LSR}),  the solution path of the TOGP was obtained as of type LSR, presented in Fig.~\ref{Fig:case3a}.  Then,  the reverse coordinate transformation in  \cite{Yuan:2020,Bakolas:2013}  was used,  leading the path of LSR into the shortest path in constant drift field,  presented in Fig.~\ref{Fig:case3a1}.
The optimal control strategy is reported in Fig.~\ref{Fig:case3c}.
It is seen from Fig.~\ref{Fig:case3a1} that the solution path in constant drift field is not the concatenation of circular arcs and straight line segments.





\subsection{Case C: Time-Optimal Guidance with a Moving Target in Constant Drift Field}

For case C, we consider a more complex problem that is to find time-optimal fixed-impact-angle guidance law for intercepting a moving target with the effect of a constant drift. Assume the velocity of the constant drift is $\boldsymbol{w}\in \mathbb{R}^2$, and let $\boldsymbol{v}\in \mathbb{R}^2$ denote the velocity of the moving target.  We consider $\boldsymbol{w}=(0.5, 0.5)$ m/s and ${\boldsymbol{v}}=(0.6368, 0.8759)$ m/s.  The initial position of the target is $(\bar{x}_0,\bar{y}_0)=(-5.5535, -0.6391)$ m, and the final heading angle is  set as $\theta_{P_f}=2.7925$ rad.

\begin{figure}[!htp]
\centering
\begin{subfigure}[t]{0.48\textwidth}
\centering
\includegraphics[width = \textwidth,trim=2cm 0cm 2cm 0cm]{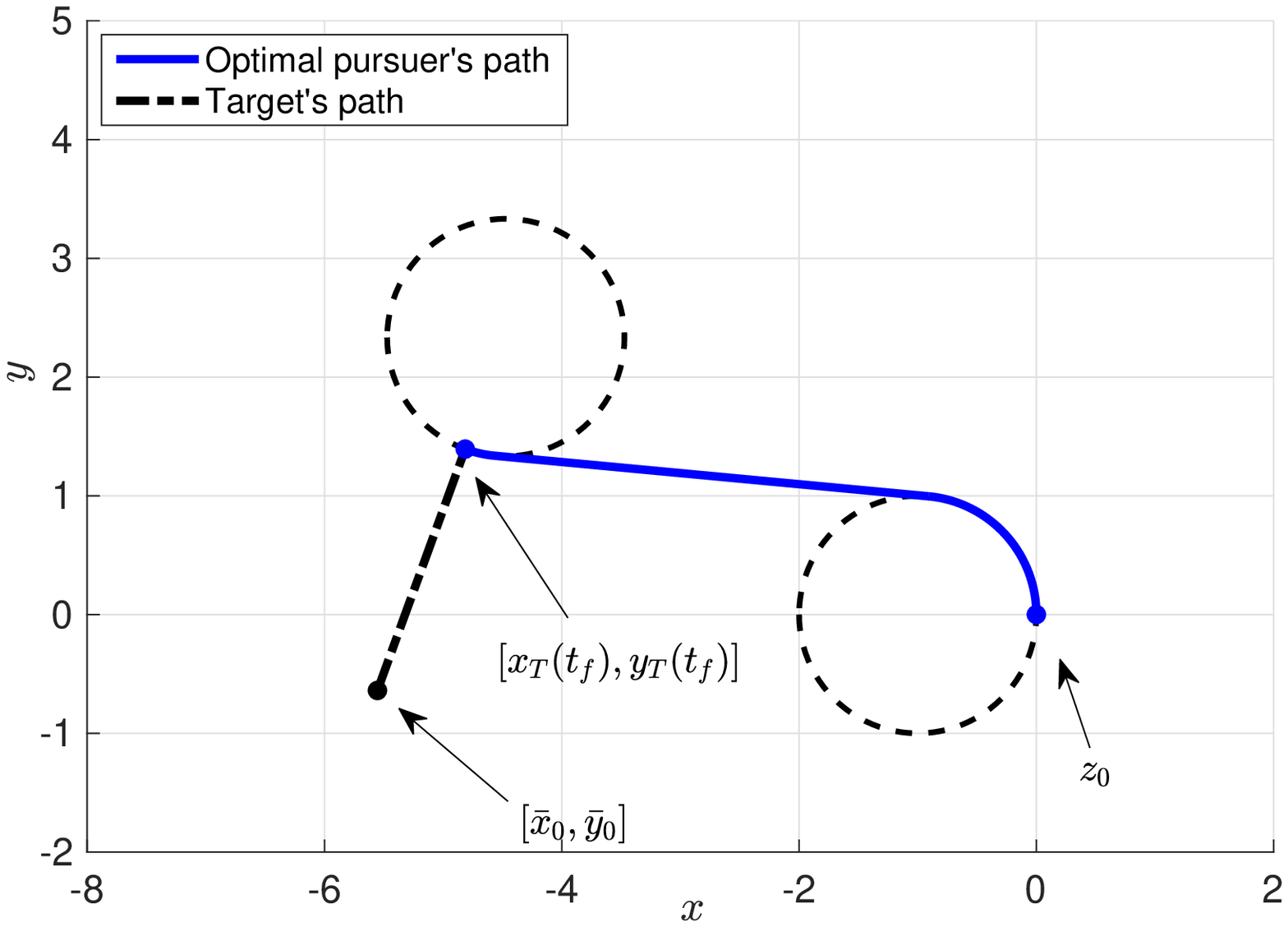}
\caption{Target's velocity is $\boldsymbol{v}-\boldsymbol{w}$ without constant drift}
\label{Fig:case4a}
\end{subfigure}
\begin{subfigure}[t]{0.48\textwidth}
\centering
\includegraphics[width = \textwidth,trim=2cm 0cm 2cm 0cm]{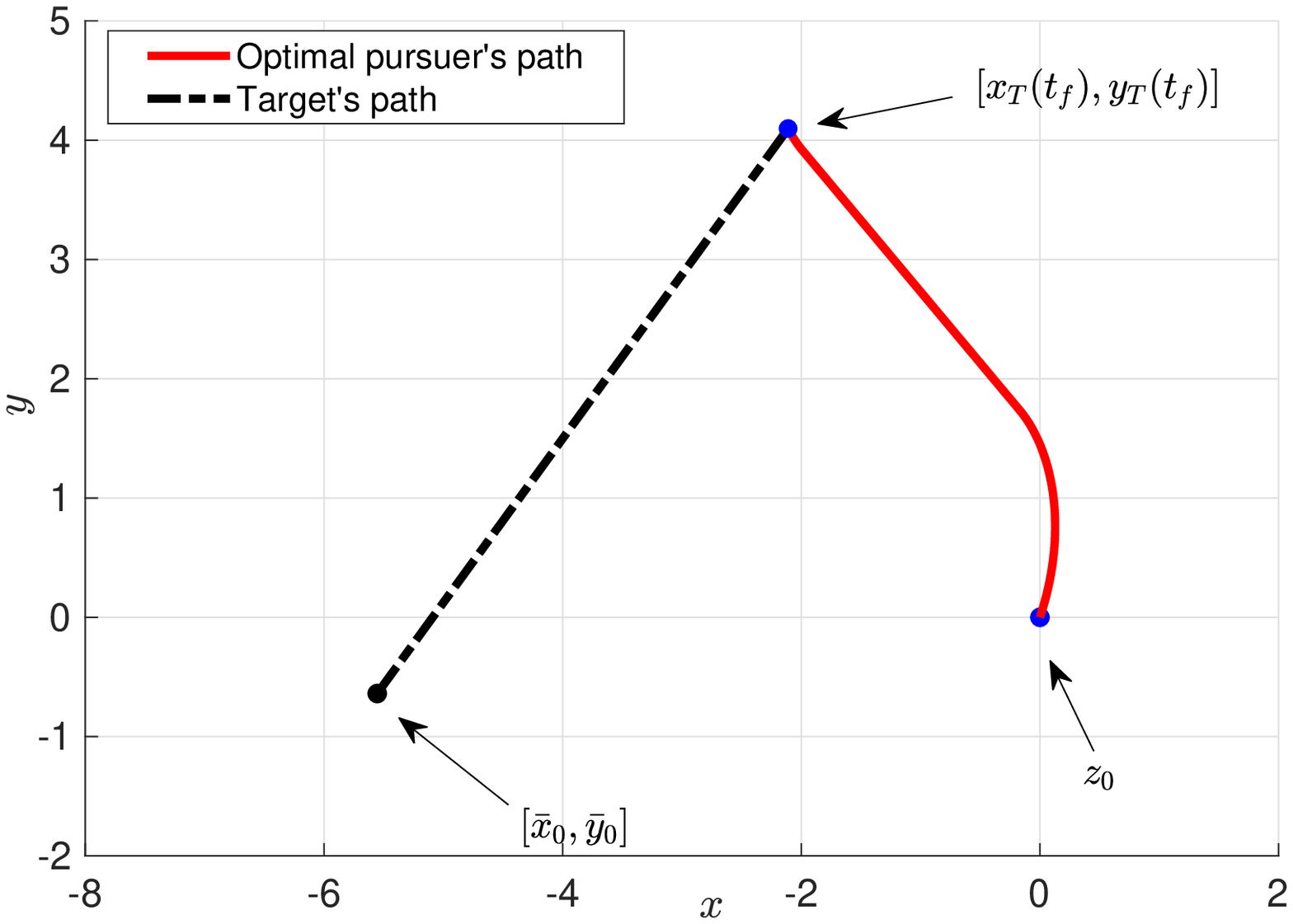}
\caption{Target's velocity is $\boldsymbol{v}$ with constant drift vector $\boldsymbol{w}$}
\label{Fig:case4a1}
\end{subfigure}
\caption{Case C: The solution path of the OCP and the corresponding path in constant drift.}
\label{Fig:CaseC:solution}
\end{figure}

\begin{figure}[t]
\centering
\includegraphics[width = 0.6\textwidth,trim=2cm 0cm 2cm 0cm]{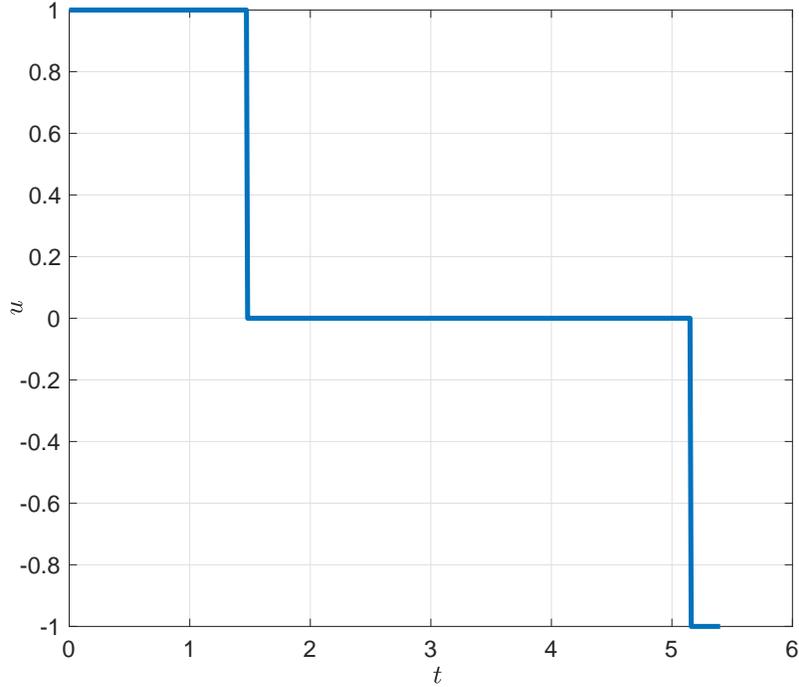}
\caption{Case C: The optimal control profile against time.}
\label{Fig:case4c}
\end{figure}

To the authors' best knowledge, this problem has not be addressed in the literature. 
By the coordinate transformation in \cite{Yuan:2020,Bakolas:2013}, this problem is equivalent to the TOGP for intercepting a moving target with velocity of $-\boldsymbol{w}+\boldsymbol{v}$. 
As a result, we are able to employ the numerical procedure in Section \ref{SE:Analytical} to find the solution of the TOGP, presented in Fig.~\ref{Fig:case4a}. 
By the reverse of the coordinate transformation \cite{Yuan:2020,Bakolas:2013}, we immediately have the time-optimal path of the original guidance problem in constant drift of $\boldsymbol{w}$,  presented in Fig.~\ref{Fig:case4a1}. The optimal control strategy is reported in Fig.~\ref{Fig:case4c}.
It should be noted that the optimal control strategy for both paths in Fig.~\ref{Fig:CaseC:solution} are the same. 

According to the procedure in solving this example, it is concluded that the developments of the paper allow efficiently and robustly finding the time-optimal fixed-impact-angle guidance law for intercepting moving targets even if the motion of the pursuer is affected by a constant drift.


\section{Conclusions}
This paper is concerned with devising the time-optimal guidance law for intercepting a moving target with a fixed impact angle. The solution paths of such time-optimal guidance problems were characterized, showing that under a reasonable assumption on the distance between initial and final positions, the solution path must lie in a sufficient family of four candidates. The geometric property of each candidate was used to formulate a nonlinear equation in terms of the candidate's parameters. As a result,  computing the time-optimal guidance law was transformed to finding the zeros of those nonlinear equations.  As each nonlinear equation might have multiple zeros and existing numerical solvers could not be guaranteed to converge to the desired zero related to the optimal path,  an efficient and robust method was proposed to find all the zeros of those nonlinear equations. The optimal path could be obtained by ruling out useless zeros.  Numerical simulations showed that  the time-optimal guidance law for intercepting moving targets could be found efficiently and robustly in comparison with the existing methods in the literature. In addition, it was also shown that the developments of this paper allowed  finding time-optimal fixed-impact-angle guidance law for intercepting moving targets even if the motion of pursuer was affected by constant drift.  


\section*{Appendix}
\appendix

\section{Expressions of $\alpha$, $\gamma$, and $d$ in terms of $\beta$}\label{Appendix:A}

In the following four subsections, we shall show that for each type in C$_{\alpha}$S$_{d}^{\beta}$C$_{\gamma}$, the values of $\alpha$, $\gamma$, and $d$ are totally determined by $\beta$. Before proceeding, we denote by $\boldsymbol{c}_0^r\in \mathbb{R}^2$ and $\boldsymbol{c}_0^l\in\mathbb{R}^2$ the centers of the right and left circles tangent to initial position and initial velocity, respectively, and denote by $\boldsymbol{c}_f^r\in \mathbb{R}^2$ and $\boldsymbol{c}_f^l\mathbb{R}^2$ the centers of the right and left circles tangent to the final position and final velocity, respectively. Then, according to the geometry in Fig.~\ref{Fig:geometry}, we immediately have
\begin{align}
\boldsymbol{c}_0^{r}  {=} 
\left(
\begin{array}{c}
\rho\\
0
\end{array}
\right) \ &\text{and} \ 
\boldsymbol{c}_0^{l}  {=}  \left(
\begin{array}{c}
-\rho\\
0
\end{array}
\right) \nonumber\\
\boldsymbol{c}_f^{r}   {=}  
\left(
\begin{array}{c}
x_f + \rho \cos (\theta_{P_f} - \pi/2)\\
y_f + \rho \sin (\theta_{P_f} - \pi/2)
\end{array}
\right) 
\ &\text{and}\ 
\boldsymbol{c}_f^{l}  {=}  
\left(
\begin{array}{c}
x_f + \rho \cos (\theta_{P_f} + \pi/2)\\
y_f + \rho \sin (\theta_{P_f} + \pi/2)
\end{array}
\right)\nonumber
\end{align}
where $(x_f,y_f)\in\mathbb{R}^2$ denotes the final position of the engagement. Set $v_x = V_T cos \theta_T$ and $v_y = V_T \sin \theta_T$. It is clear that $(v_x,v_y)$ is the constant velocity of the moving target. Then, we have the expressions of $\alpha$, $\gamma$, and $d$ for each type of C$_{\alpha}$S$_{d}^{\beta}$C$_{\gamma}$.

\subsection{Expressions $\alpha$, $\gamma$, and $d$ for RSR}

In view of  the geometry in Fig.~\ref{Fig:geometry_RSR}, if the initial circular arc has a right-turning direction, we have
\begin{align}
\alpha =
\begin{cases}
\pi/2 - \beta, \ \text{if}\ \beta\in [0,\pi/2]\\
5\pi/2 - \beta, \ \text{if}\ \beta \in (\pi/2,2\pi)
\end{cases}
\label{EQ:RSR_alpha_beta}
\end{align}
If the final circular arc has a right-turning direction, we have
\begin{align}
\gamma =
 \begin{cases}
\beta - \theta_{P_f},\ \text{if}\ \theta_{P_f} \in [0,\beta]\\
\beta - \theta_{P_f} + 2\pi, \ \text{if}\ \theta_{P_f} \in (\beta,2\pi)
\end{cases}
\label{EQ:RSR_gamma_beta}
\end{align}
In view of Fig.~\ref{Fig:geometry_RSR}, we have
\begin{align}
\begin{split}
\boldsymbol{c}_f^r &= \boldsymbol{c}_0^r + d\left(
\cos \beta,
\sin \beta
\right) ^T\\
 \frac{v_y}{v_x}&= \frac{y_f - \bar{y}_0}{ x_f - \bar{x}_0}
\end{split}
\label{EQ:RSR0_d}
\end{align}
Rearranging these two equations to eliminate $x_f$ and $y_f$ leads to
\begin{align}
d =\frac{\frac{v_y}{v_x}\rho-\frac{v_y}{v_x}\rho \sin (\theta_{P_f})-\frac{v_y}{v_x}\bar{x}_0 + \bar{y}_0 - \rho \cos (\theta_{P_f})}{\sin \beta -\frac{v_y}{v_x} \cos \beta}
\label{EQ:RSR_d_beta}
\end{align}
Up to now, it has been apparent from Eqs.~(\ref{EQ:RSR_alpha_beta}--\ref{EQ:RSR_d_beta}) that, if the path is of type RSR, the values of $\alpha$, $\gamma$, and $d$ are determined by $\beta$.


\subsection{Expressions of $\alpha$, $\gamma$, and  $d$ for LSL}

According to Fig.~\ref{Fig:geometry_LSL}, if the initial circular arc has a left-turning direction, we have
\begin{align}
\alpha = \begin{cases}
 \beta - \pi/2, \ \text{if}\ \beta\in [\pi/2,2\pi)\\
\beta + 3\pi/2, \ \text{if}\ \beta \in [0,\pi/2)
\end{cases}
\label{EQ:LSL_alpha_beta}
\end{align}
If the final circular arc has a left-turning direction, we have
\begin{align}
\gamma = 
\begin{cases}
\theta_{P_f} -  \beta, \ \text{if}\ \beta\in [0,\theta_{P_f}]\\
\theta_{P_f} - \beta + 2\pi,\ \text{if}\ \beta\in (\theta_{P_f},2\pi)
\end{cases}
\label{EQ:LSL_gamma_beta}
\end{align}
It is clear from Fig.~\ref{Fig:geometry_LSL} that the following two equations hold.
\begin{align}
\begin{split}
\boldsymbol{c}_f^l &= \boldsymbol{c}_0^l + d\left(
\cos \beta,
\sin \beta
\right) ^T\\
\frac{v_y}{v_x} &= \frac{y_f - \bar{y}_0}{ x_f - \bar{x}_0} 
\end{split}
\label{EQ:LSL0_d}
\end{align}
Eliminating $x_f$ and $y_f$ from the two equations yields
\begin{align}
d =\frac{-\frac{v_y}{v_x}\rho+\frac{v_y}{v_x}\rho \sin (\theta_{P_f})-\frac{v_y}{v_x}\bar{x}_0 + \bar{y}_0 + \rho \cos (\theta_{P_f})}{\sin \beta -\frac{v_y}{v_x} \cos \beta}
\label{EQ:LSL_d_beta}
\end{align}
We can see from Eqs.~(\ref{EQ:LSL_alpha_beta}-\ref{EQ:LSL_d_beta}) that for the path of type LSL the values of $\alpha$, $\gamma$, and $d$ are determined by $\beta$.

\subsection{Expressions of $\alpha$, $\gamma$, and $d$ for RSL}

If the path is of type RSL, the expressions of $\alpha$ and $\gamma$ are given by Eq.~(\ref{EQ:RSR_alpha_beta}) and Eq.~(\ref{EQ:LSL_gamma_beta}), respectively.  According to Fig.~\ref{Fig:geometry_RSL}, we have
\begin{align}
\begin{split}
 \boldsymbol{c}_f^l &= \boldsymbol{c}_0^r + 2 \rho \left(
\cos(\beta + \pi/2),
\sin(\beta+ \pi/2)
\right)^T
+ d \left(
\cos(\beta),
\sin(\beta)\right)^T\\
{v_y}/{ v_x}&=({y_f - \bar{y}_0})/({x_f - \bar{x}_0})
\end{split}
\label{EQ:RSL0_d}
\end{align}
Eliminating $x_f$ and $y_f$ from the two equations, we have
\begin{align}
 d =\frac{\frac{v_y}{v_x}\rho - 2\frac{v_y}{v_x}  \rho \sin \beta -\frac{v_y}{v_x}\rho \cos (\theta_{P_f} + \frac{\pi}{2}) - 2\rho \cos \beta -\frac{v_y}{v_x}\bar{x}_0 + \bar{y}_0 + \rho \sin (\theta_{P_f} +\frac{\pi}{2}) }{\sin \beta -\frac{v_y}{v_x} \cos \beta}
\label{EQ:RSL_d_beta}
\end{align}
Therefore, if the path is of type RSL, we can use Eq.~(\ref{EQ:RSR_alpha_beta}), Eq.~(\ref{EQ:LSL_gamma_beta}), and Eq.~(\ref{EQ:RSL_d_beta}) to determine the values of $\alpha$, $\gamma$, and $d$, respectively, and they are apparently determined by $\beta$.

\subsection{Expressions of $\alpha$, $\gamma$, and  $d$ for LSR}

If the path is of type LSR, the expressions of $\alpha$ and $\gamma$ are given by Eq.~(\ref{EQ:LSL_alpha_beta}) and Eq.~(\ref{EQ:RSR_gamma_beta}), respectively. According to Fig.~\ref{Fig:geometry_LSR}, we have
\begin{align}
\begin{split}
\boldsymbol{c}_f^r &= \boldsymbol{c}_0^l + 2\rho 
\left(
\cos (\beta- \pi/2),
\sin(\beta - \pi/2)
\right)^T
+ d 
\left(
\cos (\beta),
\sin(\beta )
\right)^T\\
{v_y}/{v_x} &= ({y_f - \bar{y}_0})/({x_f - \bar{x}_0}) 
\end{split}
\label{EQ:LSR0_d}
\end{align}
Eliminating $x_f$ and $y_f$ from the two equations, we have
\begin{align}
 d =\frac{-\frac{v_y}{v_x}\rho + 2\frac{v_y}{v_x}  \rho \sin \beta -\frac{v_y}{v_x}\rho \cos (\theta_{P_f} - \frac{\pi}{2}) + 2\rho \cos \beta -\frac{v_y}{v_x}\bar{x}_0 + \bar{y}_0 + \rho \sin (\theta_{P_f} -\frac{\pi}{2}) }{\sin \beta -\frac{v_y}{v_x} \cos \beta}
\label{EQ:LSR_d_beta}
\end{align}
It is apparent that for the path of type LSR we can use Eq.~(\ref{EQ:LSL_alpha_beta}), Eq.~(\ref{EQ:RSR_gamma_beta}), and Eq.~(\ref{EQ:LSR_d_beta}) to determine the values of $\alpha$, $\gamma$, and $d$, respectively, and they are determined by $\beta$ explicitly.

\section{Proof of Theorem \ref{TH:1}}\label{Appendix:B}

In the following paragraphs, each statement of Theorem \ref{TH:1} will be proven independently.
\subsection{Proof of the First Statement}

Note that the time for the intercepter from the initial condition $\boldsymbol{z}_0$ to the interception point $(x_f,y_f)$ is the same as that for the target from its initial point $(\bar{x}_0,\bar{y}_0)$ to the interception point, indicating
\begin{align}
\rho (\theta_{P_0}- \beta + \beta - \theta_{P_f} + 2n \pi) + d = \frac{\sqrt{(x_f - \bar{x}_0)^2 + (y_f - \bar{y}_0)^2}}{\sqrt{v_x^2 + v_y^2}}
\label{EQ:RSR1}
\end{align}
where $n=0$ if $\theta_{P_f}<\pi/2$ and $n=1$ if $\theta_{P_f}\geq \pi/2$. The second equation of Eq.~(\ref{EQ:RSR0_d}) is equivalent to
\begin{align}
y_f=\frac{v_y}{v_x}x_f-\frac{v_y}{v_x}\bar{x}_0+\bar{y}_0
\label{EQ:RSR2}
\end{align}
Combining Eq.~(\ref{EQ:RSR2}) and Eq.~(\ref{EQ:RSR1}) yields
\begin{align}
d=\frac{x_f-\bar{x}_0}{v_x}-\rho (\theta_{P_0}-\theta_{P_f}+2n\pi)
\label{EQ:RSR3}
\end{align}
Substituting Eq.~(\ref{EQ:RSR2}) and Eq.~(\ref{EQ:RSR3}) into Eq. ~(\ref{EQ:RSR0_d}) to eliminate $x_f$ and $d$, we have
\begin{align}\label{EQ:RSR4}
a_1 + a_2 \sin \beta +a_3\cos \beta  =0
\end{align}
where
\begin{align}
\begin{split}
a_1&=(\rho - \rho \sin(\theta_{P_f})) v_y - (\frac{v_y}{v_x}\bar{x}_0-\bar{y}_0+\rho \cos(\theta_{P_f})) v_x\\
a_2&=-\rho + \rho \sin(\theta_{P_f})-(-\frac{\bar{x}_0}{v_x}-\rho (\theta_{P_0}-\theta_{P_f}+2n\pi))v_x\\
a_3&=\frac{v_y}{v_x}\bar{x}_0-\bar{y}_0+\rho \cos(\theta_{P_f})+(-\frac{\bar{x}_0}{v_x}-\rho (\theta_{P_0}-\theta_{P_f}+2n\pi)) v_y\\
\end{split}\notag
\end{align}

\subsection{Proof of the Second Statement}
Note that the time for the intercepter from its initial condition $\boldsymbol{z}_0$ to the interception point $(x_f,y_f)$ is the same as that for the target from its initial point $(\bar{x}_0,\bar{y}_0)$ to the interception point, which means
\begin{align}
\rho (\beta - \theta_{P_0}+ \theta_{P_f} - \beta + 2n \pi) + d =
 \frac{\sqrt{(x_f - \bar{x}_0)^2 + (y_f - \bar{y}_0)^2}}{\sqrt{v_x^2 + v_y^2}}
\label{EQ:LSL1}
\end{align}
where $n=0$ if $\theta_{P_f}\geq \pi/2$ and $n=1$ if $\theta_{P_f}< \pi/2$.

By the same procedure as proving the first statement, we can combine Eq.~(\ref{EQ:LSL0_d}) and Eq.~(\ref{EQ:LSL1}) to suppress $d$, $x_f$, and $y_f$, which leads to
\begin{align}
	b_1+b_2 \sin \beta+ b_3 \cos \beta=0
\end{align}
where 
\begin{align}
    \begin{split}
    b_1&=(-\rho+\rho \sin(\theta_{P_f}))v_y-(\frac{v_y}{v_x}\bar{x}_0-\bar{y}_0-\rho \cos (\theta_{P_f}))v_x\\
    b_2&=\rho-\rho \sin(\theta_{P_f})+(\frac{\bar{x_0}}{v_x} + \rho (\theta_{P_f}-\theta_{P_0}+2n\pi))v_x\\
    b_3&=(-\frac{\bar{x}_0}{v_x} - \rho (\theta_{P_f}-\theta_{P_0}+2n\pi))v_y+\frac{v_y}{v_x}\bar{x}_0-\bar{y}_0-\rho \cos (\theta_{P_f})\\
    \end{split}
	\notag
\end{align}


\subsection{Proof of the Third Statement} 

As the time for the intercepter from its initial condition $\boldsymbol{z}_0$ to the interception point $(x_f,y_f)$ is the same as that for the target from its initial point $(\bar{x}_0,\bar{y}_0)$ to the interception point, we have
\begin{align}
\rho (\theta_{P_0}- \beta + \theta_{P_f} - \beta + 2n \pi) + 
d=\frac{\sqrt{(x_f - \bar{x}_0)^2 + (y_f - \bar{y}_0)^2}}{\sqrt{v_x^2 + v_y^2}}
\label{EQ:RSL1}
\end{align}
where
\begin{align}
n=\begin{cases}
0 & \text{if}\ \theta_{P_0}>\beta \  and \ \theta_{P_f}>\beta\\
1 & \text{if}\ (\theta_{P_0}-\beta)(\theta_{P_f}-\beta)<0\\
2 & \text{if}\ \theta_{P_0}<\beta\ and\  \theta_{P_f}<\beta
\end{cases}\nonumber
\end{align}
Rewriting the second equation of Eq.~(\ref{EQ:RSL0_d}) yields
\begin{align}
y_f=\frac{v_y}{v_x} x_f -\frac{v_y}{v_x} \bar{x}_0 + \bar{y}_0
\label{EQ:RSL2}
\end{align}
Substituting Eq.~(\ref{EQ:RSL2}) into Eq.~(\ref{EQ:RSL1}) to eliminate $y_f$ leads to
\begin{align}
d=\frac{x_f-\bar{x}_0}{v_x} - \rho (\theta_{P_0}+\theta_{P_f}-2\beta+2n\pi)
\label{EQ:RSL3}
\end{align}
Substituting Eq.~(\ref{EQ:RSL2}) and Eq.~(\ref{EQ:RSL3}) into Eq.~(\ref{EQ:RSL0_d}) to eliminate $x_f$, we eventually have
\begin{align}
c_1 +c_2 \sin \beta +c_3 \cos \beta +\beta(c_4  \sin \beta +c_5  \cos \beta) =0
\label{EQ:RSL5}
\end{align}
where
\begin{align}
\begin{split}
c_1 &=2\rho+(\rho + \rho \sin(\theta_{P_f}))v_y-(- \rho \cos(\theta_{P_f}) + \frac{v_y}{v_x} \bar{x}_0 - \bar{y}_0)v_x\\
c_2 &=-2\rho v_y-(-\frac{\bar{x}_0}{v_x}-\rho (\theta_{P_0}+ \theta_{P_f} + 2n \pi))v_x - \rho  - \rho \sin(\theta_{P_f})\\
c_3 &=-2\rho v_x +(-\frac{\bar{x}_0}{v_x}-\rho (\theta_{P_0}+ \theta_{P_f} + 2n \pi))v_y   - \rho \cos(\theta_{P_f}) + \frac{v_y}{v_x} \bar{x}_0 - \bar{y}_0\\
c_4 &=2\rho v_y\\
c_5 &=-2\rho v_x
\end{split}
\notag
\end{align}

\subsection{Proof of the Fourth Statement}


As the time for the intercepter from its initial condition $\boldsymbol{z}_0$ to the interception point $(x_f,y_f)$ is the same as that for the target from its initial point $(\bar{x}_0,\bar{y}_0)$ to the interception point, we have
\begin{align}
\rho(\beta - \theta_{P_0}+ \beta -\theta_{P_f} + &2n\pi) + d =\frac{\sqrt{(x_f - \bar{x}_0)^2 + (y_f - \bar{y}_0)^2}}{\sqrt{v_x^2 + v_y^2}}
\label{EQ:LSR1}
\end{align}
where 
\begin{align}
n=\begin{cases}
0 & \text{if}\ \theta_{P_0}<\beta \  and \ \theta_{P_f}<\beta\\
1 & \text{if}\ (\theta_{P_0}-\beta)\cdot (\theta_{P_f}-\beta)<0\\
2 & \text{if}\ \theta_{P_0}>\beta\ and\  \theta_{P_f}>\beta
\end{cases}\nonumber
\end{align} 
Rewriting the second equation of Eq.~(\ref{EQ:LSR0_d}) leads to
\begin{align}
y_f=\frac{v_y}{v_x} x_f - \frac{v_y}{v_x} \bar{x}_0 + \bar{y}_0
\label{EQ:LSR2}
\end{align}
Substituting Eq.~(\ref{EQ:LSR2}) into Eq.~(\ref{EQ:LSR1}) yields
\begin{align}
d=\frac{x_f-\bar{x}_0}{v_x} + \rho (\theta_{P_0}+\theta_{P_f} - 2\beta+2n\pi)
\label{EQ:LSR3}
\end{align}
Substituting Eq.~(\ref{EQ:LSR2}) and Eq.~(\ref{EQ:LSR3}) into Eq.~(\ref{EQ:LSR0_d}) to eliminate $x_f$, we eventually have
\begin{align}
d_1 +d_2 \sin \alpha +d_3 \cos \alpha +\alpha (d_4  \sin \alpha + d_5 \cos \alpha) =0
\label{EQ:LSR5}
\end{align}
where
\begin{align}
\begin{split}
d_1 &=-2\rho + (- \rho  - \rho \sin(\theta_{P_f}))v_y - (\rho \cos(\theta_{P_f}) + \frac{v_y}{v_x} \bar{x}_0 - \bar{y}_0)v_x\\
d_2&=2\rho v_y-(-\frac{\bar{x}_0}{v_x}+\rho (\theta_{P_0}+ \theta_{P_f} + 2n \pi))v_x +\rho  + \rho \sin(\theta_{P_f})\\
d_3&=2\rho v_x+(-\frac{\bar{x}_0}{v_x}+\rho (\theta_{P_0}+ \theta_{P_f} + 2n \pi))v_y+\rho \cos(\theta_{P_f}) + \frac{v_y}{v_x} \bar{x}_0 - \bar{y}_0\\
d_4&=2\rho v_x\\
d_5&=-2\rho v_y\\
\end{split}
\notag
\end{align}

\bibliographystyle{unsrt}  
\bibliography{references}  






\end{document}